\documentstyle[12pt]{article}

\newcommand{\wt}{\widetilde}

\newcommand{\cG}{{\cal G}}

\begin{document}

\title{\bf Approximation probabilities,
the law of quasistable markets,
and phase transitions from 
the ``condensed" state}
\author{V.P.Maslov}
\date{}
\maketitle

\begin{abstract}
For common people, in contrast to brokers, 
bankers, and those who play 
on rising and falling prices of stocks,  
the stock market law is based on the simple fact 
that the depositors aim for financial profit
at any given concrete stage.  
The common depositor cannot cause
any significant variations in prices. 
This concept suggests an analogy with the quasistable physics,
i.e., thermodynamics, in the situation in which 
the temperature varies slowly along with  the 
external conditions. 
Therefore, 
in the quasistable market,
we can see phase transitions similar to those 
in the situation of the Bose-condesate 
in thermodynamics.
We stress the positive role of information 
for common depositors 
and the possibility of changing
bonds of large denomination 
into bonds of small denomination.
\end{abstract}

We consider a discrete set of random variables
taking values $x_1,\dots,x_n$ 
with pro\-ba\-bi\-li\-ties  $p_1,\dots,p_n$.

Suppose that a number~$M$ is given.
This number will be called the (mean) number of tests.
We set $P_i = Mp_i$, so that $\sum P_i =M$.

The concept stated below is developed 
from an analysis of the quasistable market 
discussed with V.~N.~Baturin and S.~G.~Lebedev,
and from an analysis of the market approximation theory
developed by B.~S.~Kashin.

First, as was already noted by the author in his previos
works~\cite{Theorver, Maslov1, Maslov2,Maslov3, Maslov4},
people distinguish the money bonds 
they have in circulation 
only by the nomial cost of bonds,
but not by their number sign. 
In other words, although the bonds are,  
in principle, different (in number signs), 
in the market and financial problems,  
it is possible to assume that 
they are indistinguishable and hence 
subject to the Bose--Einstein statistics.
Thus, from the very beginning,  
we start from a concept other than 
that from which the standard theory of probabilities
originates.

We are concerned with a stock holder,
i.e., with a common person 
who does not play on rising price
and does not risk, 
but is looking for a direct financial income.
We study a quasistatistical market,
i.e., a market that varies slowly 
and is stationary on some time interval.

{\bf Example 1.}

The depositor has two possibilities: 
1) to deposit the money in a single $(G_1=1)$ pyramid 
(say, the ``MMM'' bank); \ \ 
2) to deposit the money in one $G_2 > 1$ 
of the common banks with the same bankrate.
The depositor has~$K$ \ 1000-rouble bonds.
Obviously, 
the number of possible deposit versions is equal to 
$$
C^K_{G+K+1} = \frac {(G+K-1)!} {(G-1)!K!}.
$$
The total number of versions is 
\begin{equation}
\prod_{k_1+k_2=K} \frac{(G_1+k_1 -1)!(G_2+k_2-1)!}
{(G_1-1)!(k_1)!(G_2-1)!k_2!} = \prod^K_{k_2=1}
\frac{(G_2+k_2-1)!}{(G_2-1)! k_2!}. 
\label{1'1}
\end{equation}
We assign the value~$x_1$ to the pyramid
and the value~$x_2$ to banks.

We assume that the average income is $\sum^2_{i=1} x_i k_i$, 
which, in particular, means that $x_2 \gg x_1$, 
since the income obtained from the pyramid
is much larger than that obtained in a bank.
Moreover, precisely as in card games, 
the price of a chip can vary, 
since, in our situation, 
it is independent of the bond denomination
(it can be a 1000-rouble bond or, e.g., a \$100-bond). 
We denote the denomination of a chip
by the letter~$\beta$, i.e.,
the gain is equal to $\beta \sum^2_{i=1} x_i  k_i$. 
The important parameter~$\beta$ varies slowly and 
independently of a common depositor.

We also have the variables~$G_1$ and~$G_2$
corresponding to~$x_1$ and~$x_2$,
which can also be treated as frequency probabilities.

Hence, according to the standard theory of probabilities,  
to the values~$x_1$ and~$x_2$ 
there correspond the probabilities
$\frac{G_1}{G_1+G_2}$, \ $\frac{G_2}{G_1+G_2}$ 
and the probabilities $\frac{k_1}{K}$, \ $\frac{k_2}{K}$.

If $k_1$, $k_2$, $G_1$, $G_2$ are large, 
then we can use the Stirling formula and rewrite~(\ref{1'1}) 
in terms of these two probabilities
(see~\cite{Theorver}).

In our example, we have $G_1=1$. 
Moreover, we can have our money not only 
in 1000-rouble bonds, 
but in any possible bonds and coins 
(with accuracy up to a 1\,copeck coin). 
Hence we shall use the approximation formula for the ``number"
of versions,  
where we replace the factorials by $\Gamma$-functions. 
Moreover, it is more convenient to consider~$\ln$ 
of the number of versions (the entropy~$S$). 
In this case, the entropy can be expressed via 
the sum or, more generally, via the Stieltjes integral.

So there are three measures corresponding to 
the values~$x_1$ and~$x_2$.

The depositor aim is to obtain the maximal income
$R= \beta \sum^2_{i=1} x_ik_i$.
Hence 
if there are neither additional taxes 
nor some additional information, 
the depositor will invest all the money in the pyramid,
i.e.,
$$
\max_{k_1+k_2=K}(\beta \sum x_i k_i) = \beta x_1 k.
$$
If the depositor has some additional information 
(the entropy~$S$) that also takes into account 
the initial frequency probabilities $G_1$ and $G_2$, 
then, in the author's oppinion, 
the following {\it main law of the market\/}
takes place, and this is the law the depositor
usually obeys in the quasistatistical stock market:
$$
R =\max_{k_1+k_2=K}\{\beta \sum^2_{i=1} x_i k_i - S(k_1, k_2)\}.
$$
This means that the depositor obtains an additional income 
from the information~$S$.

{\bf Definition.}
The expression
$$
S_M = \frac 1M \sum_i ln \Gamma(P_i)
$$
will be called the {\it symbol\/} of entropy. 
As $M \to \infty and P_i \to \infty$,
it follows from the asymptotics of the $\Gamma$-function
that $S_M \to \sum p_i \ln p_i$ as $M \to \infty$,
i.e., $S_M$ tends to the usual Shannon entropy.

In a more general form, 
this entropy looks as follows.

We set
$$
S_M(Q,P) = \frac 1M \int_\Omega \ln \Gamma (\frac{dQ}{dP})dP.
$$
If $\frac{dQ}{dP} \to \infty$ as $M \to \infty$,
then we have the limit
$$
 \lim _{M \to \infty} S_M (Q;P) = H(P,Q) 
$$
where $H(P,Q)$ is the relative entropy 
or the Kullback--Leibler information.

The quasistable market law formulated above
is similar to one of the laws in ther\-mo\-dy\-na\-mics.
The latter has not been formulated clearly as a law, 
but it is constantly used in solving different physical
problems. 
More precisely, 
the energy to which this principle leads 
was called by N.~N.~Bogoliubov 
an ``energetically efficient state."

Recall that energy states or energy levels 
in quantum mechanics are eigenvalues of some 
self-adjoint operator called the energy operator.
For simplicity, we assume that this operator is 
a finite-dimensional matrix.

Any eigenvalue is characterized by two characterstics:
one is an internal characteristic, i.e., its multiplicity, 
the other is an external characteristic.
If an eigenvalue presents energy levels, 
then this external characteristic 
is equal to the number of particles at the energy level
(or on the Bohr orbit of an atom). 
To the measure~$dP$,
there corresponds the internal characteristic
i.e., the dimension of the subspace corresponding 
to this particular level 
which is divided by the dimension of the entire space.
To the measure~$dQ$, 
there corresponds the number of particles staying 
at this level.  
Here $\int_\Omega dQ=M$ is the number of all particles.

An energetically efficient state is a state 
in which all particles are at the lower level.
If the Pauli principle 
stating that more than two particles cannot be at the same level 
is taken into account, 
then the particles must occupy all lower levels.

Thus a similar ``efficiency" principle 
also takes place in the quasistable stock market, 
but ``with a converse accuracy". 
If the largest value~$x_n$ of the random variable 
corresponds to the most profitable stocks,
then the buyer will buy all of them.
But if each bank has only one stock,
then all possible largest values $x_i$
will be bought at a rapid pace.

{\bf Example 2.}

We consider the most trivial game: 
$M$ persons play the ``heads and tails'' game with a bank. 
There are two states $\pm 1$: 
``heads'' means a gain, and ``tails'' means a loss. 
{\it Suppose that the players have the right to turn over 
the coin after it falls out.} 
Then all the players who got ``tails''
change it for ``heads'' and get in the highest place 
in the Bernoulli sequence,
although the initial probabilities are the same for each player
and the probability that all the players  get ``heads'' 
is very small.

We have considered the simplest examples
of laws of the energetical and the finanicial efficiency.

We have seen that,
in addition to the probability measures~$P$ and~$Q$,
we must introduce one more measure $\mu \ll P$
such that $\int_\Omega d\mu =K$.

In the previous papers~\cite {Theorver,Pidgin},
we introduced relative entropies corresponding to 
the Bose- and Fermi-statistics. 
They can be generalized to the case 
in which the numbers of tests are,
respectively, equal to~$M$ and~$K$ as follows:
$$
S^B_{M,K} = \int_\Omega \frac {1}{M+K}  \ln \frac 1M \Gamma
(\frac{dQ}{dP} + \frac{d\mu}{dP})- \frac 1M \ln
\Gamma(\frac{dQ}{dP}+1) -\frac 1K\Gamma(\frac{d\mu}{dP})
$$
for the Bose-statistics
(averaging over the set of stocks close in denomination) 
and
$$
S^F_{M,K} = \frac 1K \ln(\frac{d\mu}{dP}+1) - \frac 1M  \ln\Gamma
(\frac{dQ}{dP}+1) - \frac {1}{K+M}\ln\Gamma(\frac{d\mu}{dP} -
\frac{dQ}{dP}+1)
$$
for the Fermi-statistics 
(if the ``rule of queue" takes place).

So the most trivial law of the stock market
in the simplest case
says that, 
independently of the initial probabilities~$P$,
one must invest all money in the affair
that is most profitable at this particular moment.
If there are~$N$ bonds
\begin{equation}
\max (\sum P_i x_i) = Nx_n, 
\label{1'2}
\end{equation}
then all~$N$ bonds must be deposited in stocks
corresponding to the value~$x_n$ of random variables.

A similar problem is known in thermodynamics, where 
the minimum of free energy is considered.
The free energy has the form of the expression 
under the symbol~$\max$, 
where $x_i$ are energy levels 
and~$P_i$ is the number of particles at the level~$x_i$. 
This law is less transparent, but more customary.
It is universally recognized by physicists 
and confirmed by numerous experiments.

The formulas given below
are also new in thermodynamics,
but we present them for the market 
(i.e., we consider $\max$, but not~$\min$). 
It follows from the above that these formulas can be trivially
written in the language of thermodynamics.

The solution of the equation
$$
\max(\beta \sum P_i x_i - S_M)
$$
for~$P_i$ can be found from the implicit equation 
$$
\beta x_i = \int^1_0 \frac{1-z^{P_i}}{1-z} dz,
$$
which follows from the well-known formula 
for the logarithmic derivative of the $\Gamma$-function. 
Under the condition that $\sum P_i= M$,
for large~$ P_i$,
this distribution coincides with the Gibbs distribution.
We note that the physicists define the Gibbs distribution
for integer~$P_i$, but in the final formula, they obtain
noninteger~$P_i$.

In the case of a market, 
the distribution for Bose-statistics
follows from the equation
$$
\beta x_i=\int_0^1 \frac{z^{G_i+P_i-1}-z^{P_i}}{1-z} dz,
$$
and for the Fermi-statistics $P_i$, 
it follows from the equation
$$
\beta x_i=\int_0^1 \frac{z^{G_i-P_i}-z^{P_i}}{1-z} dz.
$$
In the thermodynamical case, 
in these formulas, 
$G_i$ is the multiplicity of the energy level 
of a single particle $x_i$, 
and $\beta=1/\Theta$, where $\Theta$ is the inverse temperature.
If $P_i\gg1$ and $G_i\gg 1$, 
then the solutions of these equations 
in the boson case have the form
$$
P_i\approx\frac{G_i}{e^{\beta x_i}-1},
$$
and in the fermion case, they have the form
$$
P_i\approx\frac{G_i}{e^{\beta x_i}+1}.
$$
In the thermodynamical case, the distributions 
are determined from the minimum of the expression
$$
F=\sum x_i P_i - \Theta S,
$$
which is called the free energy.
Here $S$ is the entropy. In the boson case, 
the entropy has the form
$$
S=\sum \ln\left(\frac{\Gamma(G_i+P_i)}{\Gamma(G_i)\Gamma(P_i+1)}\right),
$$
and in the fermion case, it has the form 
$$
S=\sum \ln\left(\frac{\Gamma(G_i+1)}{\Gamma(G_i-P_i+1)\Gamma(P_i+1)}\right).
$$

We shall consider the following financial model.
A depositor has some money, say, $N$,
which he can put either in the ``MMM"  bank
or in $G$ equal ``strong" banks. 
We assume that the ``MMM" bank 
gives the income $\beta\lambda_1$ per unit deposit,
while ``strong" banks give the income $\beta\lambda_2$ 
per unit deposit, where $\lambda_2<\lambda_1$ and $\beta$ 
is a positive parameter describing variations in the bankrate. 
Further, we assume that the deposit to ``strong" 
banks is equal to $k$, and the deposit to the ``MMM" bank is,
respectively, equal to $N-k$.
Then the income received by the depositor is equal to 
\begin{equation}
E(k,\beta)=\beta\lambda_1N-\beta(\lambda_1-\lambda_2)k. 
\label{1}
\end{equation}
Obviously, if there are no additional sources of income, 
it is more profitable to put all the money in the 
``MMM" bank. 
Deposits to ``strong" banks can be done in many ways.
We assume that this is related to some additional infromation
obtained by the depositor.
We also assume that any additional information
gives some additional income equal to the logarithm 
of the information amount.
Next, we consider the case in which the information amount
is determined by means of the Boltzmann statistics formula,
but here, in view of Example~1, 
we consider the case in which 
the information amount is determined by the Bose--Einstein
statistics:
\begin{equation}
\cG(k)=\frac{\Gamma(k+G)}{\Gamma(k+1)\Gamma(G)}, 
\label{2}
\end{equation}
where $\Gamma(x)$ is the Euler gamma function.
In this case, the income given by ``strong" banks 
for the deposit $k$ is equal to 
\begin{equation}
F(k,\beta)=E(k)+\ln(\cG(k))=
\beta\lambda_1N-\beta(\lambda_1-\lambda_2)k+\ln\left(
\frac{\Gamma(k+G)}{\Gamma(k+1)\Gamma(G)}\right).
\label{4}
\end{equation}
Now we study the problem of how to obtain the maximal income. 
Obviously, for this, it is necessary to find the maximum 
of function~(\ref{4}) on the interval $k\in[0,N]$.

There exists a critical value $\beta_c$ of the parameter $\beta$. 
If $\beta<\beta_c$, then function~(\ref{4}) attains its maximum
at $k=N$. This means that if all banks give low incomes,
then it is most profitable to deposit to ``strong" banks. 
The critical value is given by the formula
\begin{equation}
\beta_c=\frac{\Psi(G+N)-\Psi(N+1)}{\lambda_1-\lambda_2}, 
\label{5}
\end{equation}
where $\Psi(x)=\Gamma'(x)/\Gamma(x)$ is the derivative 
of the logarithm of the gamma function.
There exists one more critical value $\beta_0>\beta_c$
of the parameter $\beta$ such that, for $\beta>\beta_0$, 
function~(\ref{4}) has maximum at $k=0$. 
This means that, for high incomes, 
it is most profitable to give all the money to the ``MMM" bank.
This critical value is given by the formula
\begin{equation}
\beta_0=\frac{\Psi(G)-\Psi(1)}{\lambda_1-\lambda_2}. 
\label{0}
\end{equation}
If the parameter $\beta$ lies in the interval
$[\beta_c,\beta_0]$, then the maximum point of function~(\ref{4}), 
$k(\beta)$, is determined by the equation
\begin{equation}
\Psi(G+k(\beta))-\Psi(k(\beta)+1)=\beta(\lambda_1-\lambda_2). 
\label{6}
\end{equation}
By the properties of the function $\Psi(x)$, 
it is easy to see that the solution $k(\beta)$ is unique 
and is a decreasing function of~$\beta$. 

We consider expressions~(\ref{5}) and~(\ref{6})
in the limit as $N\to\infty$. We also assume that~$G$
depends on~$N$ so that the condition 
\begin{equation}
\lim_{N\to\infty}\frac{G}{N}=g>0
\label{7}
\end{equation}
is satisfied. We also take into account that 
the logarithmic derivative of the gamma function 
satisfies the relation
\begin{equation}
\Psi(G+k)-\Psi(k+1)=\int_0^1 dt\frac{t^{k}-t^{G+k-1}}{1-t}.
\label{8}
\end{equation}
Starting from~(\ref{8}), we see that~(\ref{5}) implies
\begin{equation}
\lim_{N\to\infty}\beta_c=\frac{\ln(1+g)}{\lambda_1-\lambda_2}. 
\label{9}
\end{equation}
It also follows from~(\ref{6}) that in the limit as $N\to\infty$,
under the condition~(\ref{7}), the function $k(\beta)$,
where $\beta_c<\beta$, has the form
\begin{equation}
k(\beta)=N\frac{g}{\exp(\beta(\lambda_1-\lambda_2))-1}+\mbox{O}(1). 
\label{10}
\end{equation}
In this limit case, we consider the values $\beta(m)$
for which the deposit of $m=O(N^\delta)$, $\delta<1$, 
to the ``MMM" bank gives the largest income. 
Substituting $k(\beta)=N-m$ into~(\ref{6}),
we obtain the expression
\begin{equation}
\beta(m)=\frac{\Psi(G+N-m)-\Psi(N-m+1)}{\lambda_1-\lambda_2}. 
\label{11}
\end{equation}
In the limit as $N\to\infty$ and under condition~(\ref{7}),
this expression takes the form
\begin{equation}
\beta(m)=\frac{\ln(1+g)}{\lambda_1-\lambda_2}+\mbox{O}(\frac{m}{N}). 
\label{12}
\end{equation}
Obviously, formula~(\ref{12}) implies that if
the information amout is given 
by the boson formula~(\ref{2}),
then $\beta(m)$ are close to $\beta_c$.

We show what is typical of the case in which 
the expression for the information amount 
has the Boltzmann form 
(i.e., we distinguish the bonds with different numbers
or the money is deposited by different persons).
Namely, 
\begin{equation}
\cG(k)=\frac{\Gamma(N+1)G^k}{\Gamma(k+1)\Gamma(N-k+1)}. 
\label{13}
\end{equation}
The income received in ``strong" banks $k$
takes the form 
\begin{equation}
F(k,\beta)=\beta\lambda_1N-\beta(\lambda_1-\lambda_2)k
+\ln\left(\frac{\Gamma(N+1)G^k}{\Gamma(k+1)\Gamma(N-k+1)}\right). 
\label{14}
\end{equation}
As above, we consider the problem of obtaining the maximal
income. The value of the parameter $\beta$ for which the income
is maximal, provided that the deposit to the ``MMM" bank is equal
to~$m$, is given by the formula
\begin{equation}
\beta(m)=\frac{\ln(G)+\Psi(m+1)-\Psi(N-m+1)}{\lambda_1-\lambda_2}. 
\label{15}
\end{equation}
It follows from~(\ref{15}) that, in the Boltzmann case,
there also exists a critical value $\beta_0$ such that, 
for $\beta_0\leq\beta$, 
function~(\ref{14}) attains its maximum at $k=0$:
\begin{equation}
\beta_0=\frac{\ln(G)+\Psi(N+1)-\Psi(1)}{\lambda_1-\lambda_2}. 
\label{16}
\end{equation} 
Next, since the right-hand side of~(\ref{15}) is a decreasing
function of the variable~$m$, it is obvious that the critical
value $\beta_c$, 
for which it is not profitable to have a deposit 
in the ``MMM" bank, 
exists only if the following condition is satisfied:
\begin{equation}
\ln(G)+\Psi(1)-\Psi(N+1)>0. 
\label{17}
\end{equation}
If this inequality holds, then $\beta_c$ has the form
\begin{equation}
\beta_c=\frac{\ln(G)+\Psi(1)-\Psi(N+1)}{\lambda_1-\lambda_2}.
\label{18}
\end{equation}
But if inequality~(\ref{17}) does not hold, 
then, to obtain the maximal income, 
with decreasing income given by the banks,
the deposit to the ``MMM" bank must be decreased, 
but not precisely to zero. 
As $\beta\to 0$, 
the quantity $m(\beta)$ tends to $m_{0}$, 
which is determined by the equation
\begin{equation}
\ln(G)+\Psi(m_{0}+1)-\Psi(N-m_{0}+1)=0. 
\label{19}
\end{equation}
We consider the limit as $N\to\infty$, 
assuming that condition~(\ref{7}) is satisfied. 
We have the following asymptotic relation:
\begin{equation}
\Psi(N+1)-\Psi(1)=\ln(N)+C+\mbox{0}(1), 
\label{20}
\end{equation}
where $C$ is the Euler constant. 
Therefore, inequality~(\ref{17}) has the following limit form:
\begin{equation}
g>e^C.
\label{21}
\end{equation}
For the critical value of the parameter $\beta_c$, 
provided that condition~(\ref{21}) is satisfied, 
we have
\begin{equation}
\lim_{N\to\infty}\beta_c=\frac{\ln(g)-C}{\lambda_1-\lambda_2}. 
\label{22}
\end{equation}
Moreover, obviously, if~(\ref{21}) holds, 
then, for $\beta$ at which it is profitable to have 
the deposit $m=\mbox{O}(1)$ in the ``MMM" bank,
we obtain the following asymptotics from~(\ref{15}):
\begin{equation}
\beta(m)=\frac{\ln(g)+\Psi(m+1)}{\lambda_1-\lambda_2}+\mbox{o}(1). 
\label{23}
\end{equation}
In deriving this formula, we took into account that 
$\Psi(1)=-C$. 
It follows from~(\ref{23}) that 
if the information amount is given by the Boltzmann
formula~(\ref{13}), then, in contrast to the Bose case, 
$\beta(m)$ essentially depends on~$m$.

If inequality~(\ref{21}) does not hold, 
then it follows from~(\ref{19}) that
\begin{equation}
\lim_{N\to\infty} m_{0}=\wt{m}_{0}=\mbox{O}(1)
\label{24}
\end{equation}
where $\wt{m}_{0}$ is a solution of the equation
\begin{equation}
\ln(g)+\Psi(\wt{m}_{0}+1)=0. 
\label{25}
\end{equation}
It follows from~(\ref{24}) that, 
although for a small income
it is not profitable to take away the total deposit 
from the ``MMM" bank, 
however,  
the deposit $\mbox{O}(1)$ in this bank must be small 
as compared with the total sum~$N$.

If now we consider a depositor to the pyramid and to banks, 
i.e., we assume that the bonds of the same denomination, 
but with different numbers, are ``identical",
then the phase transition, 
related to the disappearance of condensate, 
means the following for this depositor.
If $\beta$ (the price) decreases, 
then, starting from some $\beta_0$ for $\beta < \beta_0$,
the stocks cease to be sold and bought,
although it seems that it were more profitable to sell
them at any price, 
and thus somebody could speculate in selling
the stocks so that their price come to zero.
Nevertheless, in practice, this paradoxical 
fact is observed and is described in the literature.

\end{document}